\documentclass[a4paper]{article}
\usepackage{amssymb,amsmath} 

\font\tencmmib=cmmib10
\skewchar\tencmmib='177
\font\sevencmmib=cmmib7
\skewchar\sevencmmib='177
\font\fivecmmib=cmmib5
\skewchar\fivecmmib='177
\newfam\cmmibfam
\textfont\cmmibfam=\tencmmib
\scriptfont\cmmibfam=\sevencmmib
\scriptscriptfont\cmmibfam=\fivecmmib
\font\tencmbsy=cmbsy10
\skewchar\tencmbsy='60
\font\sevencmbsy=cmbsy7
\skewchar\sevencmbsy='60
\font\fivecmbsy=cmbsy5
\skewchar\fivecmbsy='60
\newfam\cmbsyfam
\textfont\cmbsyfam=\tencmbsy
\scriptfont\cmbsyfam=\sevencmbsy
\scriptscriptfont\cmbsyfam=\fivecmbsy

\begin{document}
\newtheorem{theo}{Theorem}[section]
\newtheorem{lemme}[theo]{Lemma}
\newtheorem{cor}[theo]{Corollary}
\newtheorem{defi}[theo]{Definition}
\newtheorem{problem}[theo]{Problem}
\newtheorem{prop}[theo]{Proposition}
\newtheorem{assu}[theo]{Assumption}
\newtheorem{nontheo}[theo]{Conjectured theorem}
\newcommand{\beq}{\begin{eqnarray}}
\newcommand{\enq}{\end{eqnarray}}
\newcommand{\be}{\begin{eqnarray*}}
\newcommand{\en}{\end{eqnarray*}}
\newcommand{\Td}{\mathbb T^d}
\newcommand{\T}{\mathbb T}
\newcommand{\R}{\mathbb R}
\newcommand{\N}{\mathbb N}
\newcommand{\Rd}{\mathbb R^d}
\newcommand{\Zd}{\mathbb Z^d}
\newcommand{\Linf}{L^{\infty}}
\newcommand{\dt}{\partial_t}
\newcommand{\Dt}{\frac{d}{dt}}
\newcommand{\demi}{\frac{1}{2}}
\newcommand{\ep}{^{\epsilon}}
\newcommand{\epu}{_{\epsilon}}
\newcommand{\Dtt}{\frac{d^2}{dt^2}}
\newcommand{\vf}{\varphi}
\newcommand{\bfi}{{\mathbf \Phi}}
\newcommand{\bpsi}{{\mathbf \Psi}}
\newcommand{\bm}{{\mathbf m}}
\newcommand{\NN}{\mathbb N}
\newcommand{\RR}{\mathbb R}
\newcommand{\dx}{\partial_x}
\newcommand{\vp}{v^{\perp}}
\newcommand{\E}{{\mathbf E}}
\newcommand{\Sz}{{\mathcal{S}}}
\newcommand{\ds}{\displaystyle}
\newcommand{\fe}{f_\epsilon}
\let\cal=\mathcal

\bibliographystyle{plain}
\begin{center}
\textbf{ELECTRIC TURBULENCE IN A PLASMA}\\
\textbf{SUBJECT TO A STRONG MAGNETIC FIELD}\\
\normalsize
\end{center}
\vspace{0.3cm}

\begin{center}
G. Loeper\footnotemark[1]\footnotemark[2]
A. Vasseur\footnotemark[1]\footnotemark[2]
\end{center}
\footnotetext[1]{Laboratoire J.A.Dieudonn\'e, Universit\'e de
Nice-Sophia-Antipolis, Parc Valrose, 06108 NICE Cedex 2.\\
Supported by 
}
\footnotetext[2]{Work supported by the European Atomic Energy Community EURATOM
in the Research Laboratory Agremented by CEA $\textrm{n}^o$ 01-24 UMR $\textrm{n}^o$ 
6621 CNRS-Universite de Nice.}
\vspace{0.3cm}
\begin{abstract}
We consider in this paper a plasma subject to a strong deterministic
magnetic field and we investigate the effect on this plasma
of a stochastic electric field.
We show that the limit behavior, which corresponds to the transfer of energy
from the electric wave to the particles (Landau phenomena), is described
by a Spherical  Harmonics Expansion (SHE) model.
\end{abstract}
\vspace{0.3cm}

\section{Introduction}

This paper is concerned with the effect of a stochastic electric field on
a plasma subject to a strong magnetic field. This is motivated by the study of
the electric turbulence in a fusion machine as a Tokamak.
Tokamaks are used to confine high energy plasmas in order to
obtain the conditions needed for nuclear fusion reactions to take place.
The plasma evolves in a
toroidal reactor and is confined in the heart of the torus by the the mean of a strong 
magnetic field.  A classical approximation is  to suppose  the ions to be at rest.
Then only the electrons are moving.
Another classical approximation argument  in this type of study is the
following:
we are here interested only in interactions of particles over short
distances of the order of the Larmor radius, moreover we suppose that at
this scale
the curvature of the magnetic field-lines can be neglected and  that
the plasma can be considered to be homogeneous along these field-lines.
Thus we can restrict ourselves to a bidimensional problem.
In this approach, the Vlasov equation describing the evolution of the
repartition function $f$ of the electrons is:
\begin{equation}
m\left(\frac{\partial f}{\partial t}+v\cdot \nabla_x f\right)
+q\left( B v^{\perp} + \nabla V^{\mathrm{turb}}(t,x)  \right)\cdot \nabla_{v}f=0,\label{principale1}
\end{equation}
$m$ stands for the electron's mass, $q$ its electric charge,  $f$ the distribution function
on $(t,x,v)\,\in \,\R^+\times \R^2\times\R^2$, with $t$ the time variable, $x$ the space variable and
$v$ the velocity variable. $B$ is the (constant) norm of the transverse magnetic field,
 $\nabla V^{\mathrm{turb}}$ is the turbulent electric field,
$v^{\perp}$ is the velocity vector after a rotation of  $\pi/2$.
We denote
\[
1/\epsilon=\frac{qB}{m}
\]
the cyclotronic frequency, and we want to study the effect of
 $\nabla V^{\mathrm{turb}}$
in the limit $\epsilon$ going to zero. In the deterministic case, the limits
of related problems have been studied by several authors.
In the case of the Vlasov-Poisson system (when the electric field is coupled
with the density $\int f\,dv$), the limit has been studied, even in the 3D framework,
by Fr\'enod and Sonnendrucker \cite{FS2}. In the 2D framework, using a slow
time scale adapted to the problem, the convergence of the averaged motion,
$\int_{\R^2} f(t,x,v)\,dv $ to the 2D Euler system of equations has been
 performed simultaneously by Brenier \cite{B}, Golse and Saint-Raymond
\cite{GSR} and Fr\'enod and Sonnendrucker \cite{FS}. A general result in 3D
taking into account the two effects has been performed by Saint-Raymond in
\cite{SR}.

In our case we neglect the Poisson non linear effect, concentrating on the
stochastic behavior of the equation.
Hamiltonian chaos method suggest that the modes of the turbulent electric field
 interacting with the
electrons
are those having a frequency of $\omega_n=2\pi n\epsilon$ with $n$ an integer.
This is roughly speaking the Landau resonance. Then the quasi-linear theory,
(see Garbet \cite{ref}) predicts a diffusive behavior with respect to the velocity variable.
The diffusion coefficient obtained by this method being constant, it can not take into account
the abnormal diffusion phenomena.
In this paper we are interested in turbulent electric fields
whose spectrum is spread around the Landau frequency and whose spatial
fluctuations
are of the scale of the Larmor radius (of order $\epsilon$).
We will show that the limit system is then governed by the following equation:
\begin{equation}\label{leq}
\dt \rho - \partial_e(a(e)\partial_e\rho)=0
\end{equation}
where $e=|v|^2/2$, $\rho$ is the average of $f$ over a sphere $|v|^2=2e$
and the diffusion parameter
$a(e)$ is an explicit function of the correlation of
$V^{\mathrm{turb}}$, the turbulent electric potential, and of the energy,
thus allowing abnormal diffusion.
This diffusion parameter is undimensionally defined by
(\ref{coeff_diffusion}).
This equation is similar to the so-called Spherical Harmonics Expansion (SHE)
model in high field limit modeling microelectronics semiconductor devices
(see P.Degond \cite{D} or Ben Abdallah, Degond, Markowich and Schmeiser
\cite{BADMS}).
It describes the Landau phenomena of transfer of energy from the electric wave
to the particles.
This work uses the techniques introduced by Poupaud and Vasseur \cite{PV}
to derive diffusive equation from transport in random media. This method works
directly on the equation and, for this reason, is different from the method used in 
previous works (see Kesten and Papanicolaou \cite{KP}, \cite{KP79} and 
Fannjiang, Ryzhik and Papanicolaou \cite{FRP}).
The paper is organized as follows: the precise result is stated in
Section \ref{section_results}. In Section \ref{calcul} we show how we
can compute explicitly the diffusion coefficients. Finally we give the proof
of the theorem in Section \ref{section_proof}.

\section{Results}\label{section_results}
In the remainder of the paper we fix $n$ and we denote
\begin{equation}\label{defEturb}
\nabla V^\epsilon(t,x)=\sqrt{\epsilon}\nabla V^{\mathrm{turb}}(2\pi n\epsilon t,\epsilon x)
\end{equation}
for the stochastic potential.
Equation (\ref{principale1}) takes then the following undimensional form:
\begin{eqnarray}
\frac{\partial f_\epsilon}{\partial t}+v\cdot \nabla_x f_\epsilon
+\left( \frac{v^{\perp}}{\epsilon} + \frac{1}{\sqrt{\epsilon}}\nabla V^{\epsilon}(\frac{t}{2\pi
n\epsilon},\frac{x}{\epsilon}) \right)\cdot \nabla_{v}f_\epsilon=0\label{principale}.
\end{eqnarray}
we denote $\E  $ the expectation value of any variable  and make the
following assumptions on the electrostatic potential:
\begin{eqnarray*}
&(H1)&  V^\epsilon \in  L^\infty(\RR^+;W^{3,\infty}(\RR^2)) \mbox{
and } N(\epsilon):= \E\left(
\left(\|V^\epsilon\|_{L^\infty(W^{3,\infty})}\right)^3\right) \; \; <\, \infty,\\
&(H2)&  \E V^\epsilon(t,x)  =0, \mbox{ for all  }t\in \RR^+,\;x\in \RR^2,\\
&(H3)&  V^\epsilon(t,x), \;V^\epsilon(s,y) \; \mbox {are
uncorrelated as soon as } |t-s|\geq 1,
\\
&(H4)&  \E( V^\epsilon(t,x)  V^\epsilon(s,y)) = A(t-s,x-y) +
 g^\epsilon(t,s,x,y),
\\
&&\mbox{ with: }
\\
&&\frac{\partial^{|\alpha|}}{\partial x^\alpha}A \in
L^\infty(\RR\times\RR^2) , \mbox{ for
} \alpha \in \NN^2,\; |\alpha| \leq 3,
\\
&&\| \ |\nabla^{2}_{x,y}g^\epsilon| +
|\nabla^{3}_{x,y,y}g^\epsilon| \
\|_{L^\infty((\RR^+)^2\times\RR^{4}) }
\stackrel{\epsilon\to0}{\longrightarrow} 0,
\end{eqnarray*}
where $\nabla^{2}_{x,y}g^\epsilon$ is the matrix
$\left(\partial_{x_i}\partial_{y_j}g^\epsilon\right)_{i,j}$ and
$\nabla^{3}_{x,y,y}g^\epsilon$ is
$\left(\partial_{x_i}\partial_{y_j}\partial_{y_k}g^\epsilon\right)_{i,j,k}$.

Those assumptions are the same than in \cite{PV}. 
 Hypothesis (H1) is an assumption on the regularity of
$V^\epsilon$ for $\epsilon$ fixed. Indeed the norm $N(\epsilon)$
can go to infinity when $\epsilon$ goes to 0. Hypothesis (H2)
fixes the averaged potential at 0 which is not restrictive. In
view of (\ref{defEturb}), Hypothesis (H3) determines $2\pi
n\epsilon$ as the  decorrelated lapse of time for the turbulent
electric field. Namely, it is the bigger lapse of time $t-s$ such
that the electric fields at time $t$ and at time $s$ can be
dependent on each other. Finally Hypothesis (H4), which is very
classical, can be seen as an homogeneity property which takes
place at the local scale $\epsilon$, since a quadratic quantity
which depends on four variables $(t,x,s,y)$, at the limit, depends
only on two variables $(t-s,x-y)$.

We denote  $R_s(v)$ the  rotation  of angle $s$ with center 0 of
$v$. We consider the angular average of $A$:
$$
\tilde{A}(t,x)=\frac{1}{2\pi}\int_0^{2\pi}A(t,R_\theta x)\,d\theta.
$$
We then have the following result:
\begin{theo}\label{theoreme}
Let  $V^\epsilon$ be a  stochastic potential satisfying
assumptions $(H)$ and independent of the initial data
$f^0_\epsilon\in L^2(\R^4)$. Let $a(e)$ be the function defined
by:
\begin{equation}\label{coeff_diffusion}
a(e)=\frac{1}{2\pi n^2
}\int_{0}^{+\infty}(-\partial^2_{tt}\tilde{A})(-\frac{s}{2\pi n
},2\sqrt{2e}|\sin\frac{s}{2}|)\,ds.
\end{equation}
 This function is non negative.
Assume that there is a constant  $C_0$ such that
$\|f^0_\epsilon\|_{L^2(\R^{4})}\leq C_0$ and:
\begin{equation}
\epsilon(1 + N(\epsilon)^2)\to 0.
\end{equation}
Let  $\rho_\epsilon$ be the gyro-average of $f_\epsilon$ defined by :
\begin{equation}\label{gyroaverage}
\rho_\epsilon(t,x,e)=\frac{1}{2\pi}\int_0^{2\pi}f_\epsilon(t,x,R_\theta v)\,d\theta,
\end{equation}
for every $v$ such that $|v|^2/2=e$.
Then up to extraction of a subsequence, $\E f^0_\epsilon$ converges weakly
in $L^2(\R^4)$ to a function $f^0\in L^2(\R^4)$, $\E f_\epsilon $ converges
weakly in $L^2$ to a function $f \in
L^\infty(\R^+,L^2(\R^4))$, $\E \rho_\epsilon $ converges in
$C^0([0,T],L^2(\R^{2})-w)$ for all $T>0$  toward a
function $\rho\in L^\infty(\R^+;L^2(\R^{2}\times\R^+))\cap
C^0(\R^+;L^2(\R^{2}\times\R^+)-W)$ with
$\rho(t=0,x,e)=\int_{|v|^2=2e} f^0 \,dv$. This function $\rho$
is solution to: \beq \dt \rho -
\partial_e(a(e)\partial_e\rho)=0 \ \ \ \ \ \
t>0,x\in\R^2,e\in\R^+, \label{chaleur} \enq in the distribution
sense. Finally $\E f_\epsilon $ converges weakly in $L^2(\R^+\times\R^4)$
 toward $\rho(t,x,|v|^2/2)$.
\end{theo}
Remark: depending on the regularity of the function $a(e)$ the solution of the Cauchy problem
 for equation (\ref{chaleur}) may be unique. In this case, the whole sequence $\rho_\epsilon$ converges to $\rho$ the unique solution of (\ref{chaleur}).
\section{Explicit computation of the diffusion parameter}\label{calcul}
In order to explicit the  behavior of $a(e)$ we must need the
correlation function $A$. We assume that it follows a ``Richardson-like'' law
$$
A(t,x)=f(t)|x|^\alpha.
$$
Then we have
\begin{eqnarray*}
a(e)&=&-\frac{1}{2\pi n^2}\int_{0}^{+\infty}\partial^2_{tt}\tilde{A}(\frac{-s}{2\pi n},2\sqrt{2e}|\sin\frac{s}{2}|)\,ds\\
&=&-\frac{2^{\frac{3}{2}\alpha}e^{\alpha/2}}{2\pi n^2} \int_{0}^{+\infty}f''(\frac{-s}{2\pi n})|\sin\frac{s}{2}|^\alpha\,ds\\
&=&Ke^{\alpha/2}.
\end{eqnarray*}
A necessary condition for a function of the form
$$
\rho(t,e)=\gamma(t)\rho_0(e/t^\beta)
$$
to be an auto-similar solution of  (\ref{leq}) is that
$$
\beta=\frac{2}{4-\alpha}.
$$
we then have an abnormal diffusion in
$$
e=t^{\frac{2}{4-\alpha}}.
$$
For example if $\alpha=4/3$ we find $a(e)=Ke^{2/3}$ and an abnormal diffusion in
$e=t^{3/4}$.

\section{Proof of the result}\label{section_proof}

We denote $\Sz(\R^4)$ the Schwartz space and $\Sz'(\R^4)$ its
dual. We denote $\langle\cdot;\cdot\rangle$ the duality brackets
between those two spaces. We recall that $L^2(\R^4)\subset
\Sz'(\R^4)$ and by extension we will denote
$\langle\cdot;\cdot\rangle$ as well for the scalar product on
$L^2(\R^4)$. For every linear operator $P$ on $\Sz(\R^4)$ we will
denote in the same way  $P$ its extension on $\Sz'(\R^4)$ defined
for every $\psi\in\Sz(\R^4)$ by:
\[
\langle P\psi;\eta\rangle=\langle\psi;P^*\eta\rangle, \ \ \eta\in
\Sz(\R^4).
\]
Finally we will say that $\psi_n\in\Sz'(\R^4)$ converges to
$\psi\in\Sz'(\R^4)$ in $\Sz'(\R^4)$ if for every
$\eta\in\Sz(\R^4)$,  $\langle\psi_n;\eta\rangle$ converges to
$\langle\psi;\eta\rangle$. (This is the weak convergence for
$\Sz'(\R^4)$.)

Let us rewrite equation (\ref{principale}) in the following way:
\begin{equation}\label{equationgenerale}
\left\{
\begin{array}{l}
\ds{\dt \fe+C\fe+\frac{B\fe}{\epsilon}=-\theta_t^\epsilon\fe}\\[3mm]
\ds{\fe|_{t=0}=\fe^0}
\end{array}
\right.
\end{equation}
where $C,B,\theta_t^\epsilon$ are linear operators on $\Sz(\R^4)$
defined by:
\begin{eqnarray*}
&&C=v\cdot\nabla_x\\
&&B=v^\perp\cdot\nabla_v\\
&&\theta_t^\epsilon=\frac{1}{\sqrt{\epsilon}}\nabla
V^\epsilon(\frac{t}{2\pi
n\epsilon},\frac{x}{\epsilon})\cdot\nabla_v.
\end{eqnarray*}
Notice that $C$ and $B$ are deterministic and non dependent on
$\epsilon$ nor on $t$ unlike $\theta^\epsilon_t$.

We introduce the projection operator $J$ defined on $\Sz(\R^4)$
which averages the values of the function on the spheres
$|v|^2/2=e$. Namely, for $\eta\in \Sz(\R^4)$:
\[
J\eta(x,v)=\frac{1}{2\pi}\int_0^{2\pi}\eta(x,R_\theta v)\,d\theta.
\]
We call $J$ the "gyroaverage operator". This operator is self
adjoint for the $L^2$ scalar product. Applying the projection operator $J$ on
(\ref{equationgenerale}) and taking its expectation value leads
to:
\[
\dt
\E(J\fe)+\E(JC\fe)+\frac{\E(JB\fe)}{\epsilon}=-\E(J\theta_t^\epsilon\fe).
\]
We first study some properties of the operators in order to pass
to the limit in the left hand side of this equation. Then we
investigate the limit of $\E(J\theta_t^\epsilon\fe)$ following the
procedure of \cite{PV}. Finally we derive the SHE equation giving
the explicit form of the diffusion coefficient $a(e)$.

\subsection{Properties of the operators}

We have the following properties on the operators $C$, $B$,
$\theta_t^\epsilon$ and $J$:
\begin{lemme}\label{lemmeoperator}
Operators $C$,$B$ and $\theta_t^\epsilon$ are skew adjoint for the
$L^2$ scalar product. Operators $C$, $B$, and $J$ commute with
the expectation operator $\E$. The operator $J$ is the restriction on 
$\Sz(\R^4)$ of 
the orthogonal projector of $L^2(\R^4)$ into $\mathrm{Ker}B$. In
particular:
\begin{eqnarray*}
&&\|J\|_{\mathcal{L}(L^2(\R^4))}=1,\\
&&J^2=J,\\
&&\mathrm{Ker}B=\mathrm{Im}J.
\end{eqnarray*}
In addition:
\[
JB=JCJ=0.
\]
\end{lemme}
\noindent {\bf Proof.}

\noindent
 --Operators $C$,$B$ and $\theta_t^\epsilon$
can be rewritten as $b.D$, where $D$ is a gradient operator and $b$
a regular function verifying $D\cdot b=0$. For every
 functions $\eta_1,\eta_2\in\Sz(\R^4)$ we
have:
\begin{eqnarray*}
\langle b\cdot D\eta_1;\eta_2\rangle&=&-\langle\eta_1;D(b\eta_2)\rangle\\
&=&-\langle \eta_1;b\cdot D\eta_2\rangle.
\end{eqnarray*}
Hence they are skew adjoint operators for the $L^2$ scalar
product. \vskip3mm

\noindent --The operator $J$ is clearly the $L^2$ projection on
$L^2$ functions which depends only on $|v|^2/2$ with respect to
$v$. In polar coordinates $v=(r\sin \theta,r\cos\theta)$, we have
$B=\partial/\partial\theta$. So $J$ is the projection on
$\mathrm{Ker}B$.\vskip3mm

\noindent --Operator $J$ is the projector on $\mathrm{Ker}B$,
hence $BJ=0$. Since $B$ is skew adjoint and $J$ is self adjoint,
we have $(BJ)^*=-JB=0$.\vskip3mm

\noindent --Let us fix $\eta\in\Sz(\R^4)$. We denote:
$J\eta(x,v)=\rho_\eta(x,\frac{|v|^2}{2})$. Hence
\[
JCJ\eta(x,v)=\nabla_x\cdot\left(\rho_\eta(x,\frac{|v|^2}{2})\frac{1}{2\pi}\int_0^{2\pi}
R_\theta v\,d\theta\right)=0.
\]
Finally $JCJ=0$.\vskip3mm

\noindent --Since $C$, $B$ and $J$ are linear and deterministic,
they commute with $\E$.
 \qquad$\hfill \Box$. \vskip3mm

From those properties we deduce the following proposition:
\begin{prop}\label{proposition1}
For every $\epsilon$ and every $t\in\R^+$ we have:
\[
\|\fe(t)\|_{L^2(\R^4)}= \|\fe^0\|_{L^2(\R^4)}.
\]
There exists  a function $f^0\in L^2(\R^4)$ and a function $f\in
L^\infty(\R^+;L^2(\R^4))$
such that, up
to a subsequence,  $\E\fe^0$ converges weakly in $L^2(\R^4)$ to
$f^0$, $\E\fe$ converges weakly in $L^2([0,T]\times\R^4)$ to $f$
for every $T>0$. For every $t>0$, $f$ verifies $Jf(t)=f(t)$.
The function $J\E\fe$ is solution to:
\begin{equation}\label{equationmoyenee}
\left\{
\begin{array}{l}
\ds{\dt J\E\fe+\E(J\theta_t^\epsilon\fe)=w_\epsilon}\\[3mm]
\ds{J\E\fe|_{t=0}=J\E\fe^0},
\end{array}
\right.
\end{equation}
where $w_\epsilon$ converges to 0 in $\Sz'$.
\end{prop}
\noindent {\bf Proof.} Since $C$, $B$ and $\theta_t^\epsilon$ are
skew adjoint operators,  we have:
\[
\dt\langle \fe(t);\fe(t)\rangle=0,
\]
which gives the first equality. By weak compactness there exists
two functions $f^0\in L^2(\R^4)$ and  $f\in
L^\infty(\R^+;L^2(\R^4))$ such that, up to a subsequence,
$\E\fe^0$ converges weakly in $L^2(\R^4)$ to $f^0$, $\E\fe$
converges weakly in $L^2([0,T]\times\R^4)$ to $f$ for every $T>0$.
Notice that $\epsilon \theta_t^\epsilon$ converges to 0 in
$\Sz'(\R^4)$. Multiplying Equation (\ref{equationgenerale}) by
$\epsilon$, taking its expectation value,  and letting $\epsilon$
go to 0, we find:
\[
B f(t)=0 \ \ \mathrm{on} \ \ \R^+,
\]
since $B$ and $\E$ commute. Thanks to Lemma \ref{lemmeoperator}
$f(t)\in\mathrm{Im}J$, and since $J^2=J$, we have $Jf(t)=f(t)$ for
almost every $t>0$. Since $JB=0$, applying the operator $J$ on
equation (\ref{equationgenerale}) and taking its expectation value
gives:
\[
\dt J\E\fe+\E(J\theta_t^\epsilon\fe)=w_\epsilon,
\]
with $w_\epsilon=-\E JC\fe$. This converges in $\Sz'$ to $-JC f=-JCJ
f=0$, thanks to Lemma \ref{lemmeoperator}.  \qquad$\hfill\Box$

Hence we are now concerned by the limit in $\Sz'$ of
$\E(J\theta_t^\epsilon\fe)$.

\subsection{Computation of $\E(J\theta_t^\epsilon\fe)$}

Let us denote $S^\epsilon_t \ t\in\R$ the group on $\Sz(\R^4)$
generated by the operator $C+B/\epsilon$. Namely, for every
$h\in\Sz$, $S^\epsilon_t h$ is the unique solution on $\R$ to:
\begin{equation}
\left\{
\begin{array}{l}
\ds{\dt g+Cg+\frac{Bg}{\epsilon}=0}\\
\ds{g|_{t=0}=h}.
\end{array}
\right.
\end{equation}
The operator $S^\epsilon_t$ can be explicitly  given by:
\[
S^\epsilon_th(x,v)=h(T_\epsilon(t)(x,v)),
\]
where
$$T_\epsilon(t)(x,v)=(x+\epsilon\vp -\epsilon R_{-t/\epsilon}\vp,R_{-t/\epsilon}v).$$
The function $T_\epsilon(t)(x,v)$ gives the position at $-t$ of the
particle being in $x$ with speed $v$ at time $0$ and moving at
constant speed $|v|$ on a circle of radius $\epsilon|v|$. In
particular $S^\epsilon_t$ is $2\pi\epsilon$ periodic. Notice that the
adjoint of $S^\epsilon_t$ is $S^\epsilon_{-t}$.

Following the procedure of \cite{PV}, we use a 2 times iterated
Duhamel formula.
The first iteration gives:
\be
\fe(t)=S^\epsilon_{2\pi n\epsilon } \fe(t-2\pi n\epsilon) -
\int_0^{2\pi n\epsilon} (S^\epsilon_\sigma
\theta^\epsilon_{t-\sigma}
S^\epsilon_{-\sigma})S^\epsilon_{\sigma}\fe(t-\sigma)\,d\sigma  
\en
and then we write the Duhamel formula for the $\fe(t-\sigma)$ in the integral, and this yields:
\begin{eqnarray*}
\fe(t)&=&S^\epsilon_{2\pi n\epsilon } \fe(t-2\pi n\epsilon) -
\int_0^{2\pi n\epsilon} (S^\epsilon_\sigma
\theta^\epsilon_{t-\sigma}
S^\epsilon_{-\sigma})S^\epsilon_{{4\pi n}\epsilon}\fe(t-{4\pi n}\epsilon)\,d\sigma  \\
 &&+  \int_0^{2\pi n\epsilon}\int_0^{{4\pi n}\epsilon-\sigma} (S^\epsilon_\sigma
 \theta^\epsilon_{t-\sigma}S_{-\sigma}^\epsilon)
( S_{s+\sigma}^\epsilon
 \theta^\epsilon_{t-\sigma-s}S_{-s-\sigma}^\epsilon)
S_{s+\sigma}^\epsilon \fe(t-\sigma-s) \,ds \,d\sigma.
\end{eqnarray*}

We obtain:
\begin{eqnarray}
 \nonumber\E(J\theta^\epsilon_t \fe(t))&=&J\E\left( \theta^\epsilon_t
S^\epsilon_{2\pi n\epsilon } \fe(t-2\pi n\epsilon)\right) \\
&&-
\int_0^{2\pi n\epsilon}\E\left(J\theta^\epsilon_t
(S^\epsilon_\sigma \theta^\epsilon_{t-\sigma}S^\epsilon_{-\sigma
}) S^\epsilon_{{4\pi n}\epsilon}\fe(t-{4\pi n}\epsilon) \right)
\,d\sigma +  r^\epsilon_t\label{eq:abst3}
\end{eqnarray}
with
\be
r^\epsilon_t=
\int_0^{2\pi n\epsilon}\int_0^{{4\pi n}\epsilon-\sigma} J\E\left(
\theta^\epsilon_t (S^\epsilon_\sigma
\theta^\epsilon_{t-\sigma}S^\epsilon_{-\sigma})(S^\epsilon_{s+\sigma}
\theta^\epsilon_{t-\sigma-s}S^\epsilon_{-s-\sigma})S^\epsilon_{s+\sigma}
\fe(t-\sigma-s) \right)\,ds \,d\sigma.
\en
The  function $\fe^0$ is independent of the operators
$\theta^\epsilon_t,\;t\in \RR$. In particular, in view of the
assumption (H3),  $\theta^\epsilon_t$ and $\fe(t-s)$ are independent
as soon as
$t\geq s+2\pi n\epsilon$ and $s\geq 0$.\\

Combining this fact with (H2), Equation (\ref{eq:abst3}) becomes
for $t\geq {4\pi n}\epsilon$
\begin{eqnarray} \label{eq:abst4}
\nonumber \E(J\theta^\epsilon_t\fe)&=&J\E(\theta^\epsilon_t) \E
(S^\epsilon_{2\pi n\epsilon } \fe(t-2\pi n\epsilon)) \\
&&\nonumber -
\int_0^{2\pi n\epsilon}\E(J\theta^\epsilon_t (S^\epsilon_\sigma
\theta^\epsilon_{t-\sigma}S^\epsilon_{-\sigma})) \E( S^\epsilon_{{4\pi
n}\epsilon } \fe(t-{4\pi n}\epsilon) ) \,d\sigma+ r^\epsilon_t,
\\
\E(J\theta^\epsilon_t\fe)&=& -\int_0^{2\pi n\epsilon}
\E(J\theta^\epsilon_t (S^\epsilon_\sigma
\theta^\epsilon_{t-\sigma}S^\epsilon_{-\sigma}))\E\fe(t) \,d\sigma
+ r^\epsilon_t+e^\epsilon_t,
\\
 \nonumber \mbox{with} \qquad e^\epsilon_t &=& -\int_0^{2\pi
n\epsilon}\E(J\theta^\epsilon_t (S^\epsilon_\sigma
\theta^\epsilon_{t-\sigma}S^\epsilon_{-\sigma}))(\E S^\epsilon_{{4\pi
n}\epsilon } \fe(t-{4\pi n}\epsilon) - \E \fe(t))\,d\sigma.
\end{eqnarray}

Since $S^\epsilon_t$ is $2\pi\epsilon$ periodic, $S^\epsilon_{{4\pi
n}\epsilon}\fe(t-{4\pi n}\epsilon)=\fe(t-{4\pi n}\epsilon)$. We
have:
\[
S^\epsilon_s\theta_{t-s}^\epsilon
S^\epsilon_{-s}=\frac{1}{\sqrt{\epsilon}}(S^\epsilon_sE^\epsilon(t-s))\cdot D^\epsilon_s
\]
where we denote
$$E^{\epsilon}(t,x)=\nabla V^\epsilon(\frac{t}{2\pi n\epsilon},\frac{x}{\epsilon}),$$
and we define the differential operator $D^\epsilon_s$  by
$$D_s^\epsilon=R_{-s/\epsilon}\nabla_v +\epsilon
R_{-s/\epsilon}\nabla_x^{\perp}-\epsilon\nabla_x^{\perp}.$$ 
Note that $D^\epsilon_s$ is skew adjoint.
Let us
introduce the operator $L_t^\epsilon$ on $\Sz(\R^4)$ (extended on
$\Sz'(\R^4))$ defined for every $\eta\in\Sz(\R^4)$ by:
\[
L^\epsilon_t\eta=-\int_0^{2\pi n\epsilon} \E(\theta^\epsilon_t
(S^\epsilon_\sigma
\theta^\epsilon_{t-\sigma}S^\epsilon_{-\sigma}))\eta \,d\sigma.
\]
We can gather those results in the following way:
\begin{lemme}
We have the following equality:
\[
\E(J\theta^\epsilon_t\fe)=JL_t^\epsilon\E\fe+r_t^\epsilon+e_t^\epsilon,
\]
where the operator $L^\epsilon_t$ is defined for every $\eta\in
\Sz(\R^4)$ by:
\begin{equation}
L^\epsilon_t\eta(x,v)=-\frac{1}{\epsilon}\int_0^{2\pi
n\epsilon}\nabla_v\cdot(\E(S^\epsilon_\sigma
E^\epsilon(t-\sigma)\otimes
E^\epsilon(t))\cdot D_\sigma^\epsilon\eta(x,v))\,d\sigma,
\end{equation}
and the remainders are defined by:
\begin{eqnarray*}
&&e^\epsilon_t=JL^\epsilon_t(\E\fe(t-{4\pi n}\epsilon)-\E\fe(t)),\\
&&r^\epsilon_t=\frac{1}{\epsilon\sqrt{\epsilon}}\int_0^{2\pi n\epsilon}\int_0^{{4\pi n}\epsilon-\sigma}
J\E\left( E^\epsilon(t)\cdot\nabla_v(S^\epsilon_\sigma
E^\epsilon(t-\sigma)\cdot D^\epsilon_\sigma
(S^\epsilon_{s+\sigma}E^\epsilon(t-s-\sigma)\right.\\
&&\left.\qquad\qquad \qquad\qquad\qquad\cdot D_{s+\sigma}^\epsilon(
S^\epsilon_{s+\sigma} \fe(t-\sigma-s))))\right)\,ds \,d\sigma.
\end{eqnarray*}
\end{lemme}
We can now show the following lemma:
\begin{lemme}
For every $\eta\in\Sz(\R^4)$, the remainder $r_t^\epsilon$
verifies:
\[
|\langle r_t^\epsilon;\eta\rangle|\leq
C(\eta)\sqrt{\epsilon}N(\epsilon),
\]
and $(L_t^\epsilon)^*\eta$ converges in $L^2(\R^4)$ to:
\[
\int_0^{2\pi
n}R_{-s}\nabla_v\cdot \left(\nabla^2_{xx}A(\frac{-s}{2\pi
n},v^\perp-R_{-s} v^\perp)\nabla_v \eta\right)\,ds.
\]
\end{lemme}
\noindent {\bf Proof.} We have:
\begin{eqnarray*}
(L^\epsilon_t)^*\eta &=&-\frac{1}{\epsilon}\int_0^{2\pi
n\epsilon}D^\epsilon_\sigma\cdot(\E(S_{\sigma}^\epsilon
E^\epsilon(t-\sigma)\otimes E^\epsilon(t))\cdot\nabla_v \eta)\,d\sigma\\
&=&-\int_0^{2\pi n}D^\epsilon_{\epsilon \sigma}\cdot(\E(S_{\epsilon
\sigma}^\epsilon E^\epsilon(t-\epsilon \sigma)\otimes
E^\epsilon(t))\cdot\nabla_v) \eta\,d\sigma.
\end{eqnarray*}
But thanks to the definition to $T_\epsilon(s)$, $E^\epsilon$ and
Hypothesis (H3), the term
\[
\E(S_{\epsilon \sigma}^\epsilon E^\epsilon(t-\epsilon
\sigma)\otimes E^\epsilon(t))=\E(\nabla V^\epsilon(\frac{t-\epsilon\sigma}{2\pi n
\epsilon},x/\epsilon+v^\perp-R_{-\sigma}v^\perp)\otimes
\nabla V^\epsilon(\frac{t}{2\pi n\epsilon},x/\epsilon))
\]
converges strongly to $(-\nabla^2_{xx}A)(-\sigma/(2\pi
n),v^\perp-R_{-\sigma} v^\perp)$ in $\Linf((\RR^+)^2;W^{1,\infty}(\RR^2))$. Hence thanks to the definition
of $D^\epsilon_s$, $(L^\epsilon_t)^*\eta$ converges strongly to:
\[
\int_0^{2\pi
n}R_{-s}\nabla_v\cdot \left(\nabla^2_{xx}A (\frac{-s}{2\pi
n},v^\perp-R_{-s} v^\perp)\nabla_v \eta\right)\,ds
\]
in $\Linf(\RR^+\times\RR^4)$.
We recall that we have $D^\epsilon_{\epsilon s}=R_{- s }\nabla_v
+\epsilon R_{-s}\nabla_x^{\perp}-\epsilon\nabla_x^{\perp}$ thus
\be &&\|D^\epsilon_{\epsilon s'}\Phi\|_{L^2(\R^+\times\R^4)}\leq C
\|\Phi\|_{W^{1,2}(\R^+\times\R^4)},\\
&&\|D^\epsilon_{\epsilon s'}D_{\epsilon
s}\Phi\|_{L^2(\R^+\times\R^4)}\leq C
\|\Phi\|_{W^{2,2}(\R^+\times\R^4)}, \en hence
\begin{eqnarray*}
|\langle r_t^\epsilon;\eta\rangle|
&\leq&C(n)\sqrt{\epsilon}\|\Phi\|_{W^{3,2}}\|f_\epsilon\|_{L^2}\\
&& \sup_{s,s'}\left\{ \E(|E^\epsilon(T(\epsilon
s'))|\left(|E^\epsilon(T(\epsilon s))|+|D^\epsilon_{\epsilon s'}
[E^\epsilon(T(\epsilon s))]|\right)\right.\\
&&\left.(|E^\epsilon(t,x)|+|D^\epsilon_{\epsilon
s}E^\epsilon(t,x)|+|D^\epsilon_{\epsilon
s'}E^\epsilon(t,x)|+|D^\epsilon_{\epsilon s'}D^\epsilon_{\epsilon
s}E^\epsilon(t,x)|))\right\}.
\end{eqnarray*}
We then use the following bounds: (recall that $E^\epsilon(t,x)=\nabla V^\epsilon(\frac{t}{2\pi n\epsilon},\frac{x}{\epsilon}).$)
\be
&&|E^\epsilon| \leq \|\nabla_xV^\epsilon\|_{\Linf(\R^+\times \R ^2))}\\
&&|D^\epsilon_{\epsilon s}E^\epsilon(t,x)|\leq \|\nabla^2_{xx}V^\epsilon\|_{\Linf(\R^+\times \R ^2))}\\
&&|D^\epsilon_{\epsilon s'}D^\epsilon_{\epsilon
s}E^\epsilon(t,x)|\leq
\|\nabla^3_{xxx}V^\epsilon\|_{\Linf(\R^+\times \R ^2))}\\
&&|D^\epsilon_{\epsilon s'}[E^\epsilon(T(\epsilon s))]|\leq
C\epsilon |\nabla_xE^\epsilon|(T(\epsilon s))\leq
C\|\nabla^2_{xx}V^\epsilon\|_{\Linf(\R^+\times \R ^2))}. \en Then
Hypothesis $(H1)$ ensures that
\[
|\langle r_t^\epsilon;\eta\rangle|\leq
C\sqrt{\epsilon}\|f_0^\epsilon\|_{L^2}N(\epsilon)\|\Phi\|_{W^{3,2}},
\]
which ends the proof of the lemma.
  \qquad $\hfill\Box$

We can now state the following proposition:
\begin{prop}\label{proposition2}
Assume that $\epsilon (N(\epsilon))^2$ converges to 0 when
$\epsilon$ goes to 0. Then the convergence (up to a
subsequence) of $J\E\fe$ to $f$ holds in $C^0(\R^+;L^2(\R^4)-w)$,
and $f$ is solution to:
\begin{equation}\label{equationlabonne}
\dt f+JL^0_tJf=0,
\end{equation}
where the operator $L^0_t$ is defined for every $\eta\in\Sz(\R^4)$
by:
\[
(L^0_t)^*\eta=\int_0^{2\pi
n}R_{-s}\nabla_v\cdot \left(\nabla^2_{xx}A (\frac{-s}{2\pi
n},v^\perp-R_{-s} v^\perp)\nabla_v \eta\right)\,ds.
\]
\end{prop}
\noindent {\bf Proof.} Thanks to the previous lemma, for every
test function $\eta\in \Sz(\R^4)$:
\[
|\langle
r_t^\epsilon;\eta\rangle|\stackrel{\epsilon\to0}{\longrightarrow}0,
\]
and $(L^\epsilon_t)^*\eta$ converges strongly in $L^2(\R^4)$ to $(L^0_t)^*\eta$.
But, thanks to Proposition
\ref{proposition1}, $\fe$ converges weakly to $f$ in $L^2-w.$ 
So 
$$\langle JL^\epsilon_t\E\fe;\eta\rangle
=\langle \E\fe;(L^\epsilon_t)^*J\eta \rangle$$ 
converges to:
$$\langle f;(L^0_t)^*J\eta \rangle=
\langle JL^0_tf;\eta\rangle.$$ 
The function  $\fe(t)-S_{2\pi n\epsilon}\fe(t-2\pi n\epsilon)$
converges to 0 in $L^2-w$ as well. So $e_t^\epsilon$ converges to 0 in 
$\Sz'(\R^4)$. Passing to the limit in
equation (\ref{equationmoyenee}) gives equation
(\ref{equationlabonne}). This shows that $\dt \fe$ is uniformly
bounded in time in a negative Sobolev space. Hence $\fe$ converges
to $f$ in the space of continuous function in time with values in
this Sobolev space. Finally since $\fe$ is bounded in
$L^\infty(\R^+;L^2(\R^4))$, the convergence holds in
$C^0([0,T];L^2(\R^4)-w)$ for every $T>0$. \qquad $\hfill\Box$

\subsection{Convergence to the SHE model}

Since $Jf=f$, we can introduce the gyroaverage function defined by:
\[
\rho(t,x,e)=f(t,x,v),
\]
for every $v$ such that $2e=|v|^2$. This subsection is devoted to
the proof of the following lemma:

\begin{lemme}\label{equationlimitegyromoyenee}
The function $\rho$ lies in $C^0(\R^+;L^2(\R^2\times \R^+)-w)\cap
L^\infty(\R^+;L^2(\R^2\times\R^+))$. It is solution to:
\begin{equation*}
\left\{
\begin{array}{l}
\displaystyle{\dt \rho-\partial_e(a(e)\partial_e\rho)=0}\\[3mm]
\displaystyle{\rho|_{t=0}=\frac{1}{2\pi}\int_{0}^{2\pi}  f^0
(t,x,R_\theta v) \,d\theta}
\end{array}
\right.
\end{equation*}
where the diffusion parameter is defined by:
\[
a(e)=\int_0^{2\pi n}\int_{0}^{2\pi}R_\theta
v\cdot(-\nabla^2_{xx}A)(-\frac{s}{2\pi n},R_\theta
v^\perp-R_{-s+\theta}v^\perp)\cdot R_{-s}R_\theta v\,ds\,d\theta,
\]
for every $v$ such that $e=|v|^2/2$.
\end{lemme}

\noindent {\bf Proof.} Let us first compute the operator
$JL^0_tJ$. Let $\eta_1,\eta_1$ be two test functions in
$\Sz(\R^4)$. We have:
\begin{eqnarray*}
\langle \eta_1;JL^0_tJ\eta_2\rangle&=&\langle
(L^0_t)^*J\eta_1;J\eta_2\rangle\\
&=&\int_{\R^4}\int_0^{2\pi
n}\nabla_vJ\eta_1(-\nabla^2_{xx}A)(\frac{-s}{2\pi n},v^\perp-R_{-s}
v^\perp)R_{-s}\nabla_v J\eta_2\,ds\,dx\,dv.
\end{eqnarray*}
Let us denote $\rho_{\eta_i}$ for $i=1,2$ the functions defined
by:
\[
\rho_{\eta_i}(x,\frac{|v|^2}{2})=J\eta_i(x,v).
\]
Using polar coordinates and noticing that $dv=d\theta\,de$ we
find:
\begin{eqnarray*}
\langle
\eta_1;JL^0_tJ\eta_2\rangle&=&\int_{\R^2}\int_0^\infty\partial_e\rho_{\eta_1}(x,e)
\partial_e\rho_{\eta_2}(x,e)\int_0^{2\pi n}\int_0^{2\pi}R_{\theta}
\vec{e}\cdot\\
&&\qquad\qquad (-\nabla^2_{xx}A)(\frac{-s}{2\pi
n},-R_{-s+\theta+\frac{\pi}{2}}\vec{e}+
R_{\theta+\frac{\pi}{2}}\vec{e})\cdot R_{-s+\theta}\vec{e}\,ds\,d\theta
\,de\,dx\\
&=&\int_{\R^2}\int_0^\infty\partial_e\rho_{\eta_1}(x,e)
\partial_e\rho_{\eta_2}(x,e)a(e)\,de\,dx,
\end{eqnarray*}
where $\vec{e}=(\sqrt{2e},0)$. Hence for every test function
$\rho_\eta$, let us multiply it by Equation
(\ref{equationlabonne}) and integrate with respect to $x,v$. Since
$de \,d\theta=dv$ we find:
\[
\dt \int_{\R^2}\int_0^\infty \rho(t,x,e)\rho_\eta(x,e)\,dx\,de=
\int_{\R^2}\int_0^\infty \rho(t,x,e)\partial_e
(a(e)\partial_e\rho_\eta(x,e))\,de\,dx.
\]
This, with Proposition \ref{proposition2} gives the desired
result. $\hfill\Box$
\\
\noindent Remark: we have a family of equations parametrized by
$x\in \R ^2$, and the solutions of two equations at two distinct
$x$ do not interact.

\bigskip
\noindent

\subsection{Explicit computation of the diffusion coefficient}

We derive in the following a suitable form to the diffusion
coefficient $a(e)$. We will show, in particular, that $a(e)$ is
non negative. 
From $(H4)$ the correlation function $A(t,x)$ is even with respect to $t$ and $x$. This with $(H3)$ gives:
\begin{lemme}\label{hypotheseA}
The correlation function $A$ satisfies:
\begin{eqnarray*}
&&\mathrm{Supp}A\subset[-2\pi n,2\pi n]\times\R ^2,\\
&&\nabla_x A(0,0)=0,\\
&&\partial_s A(0,0)=0.
\end{eqnarray*}
\end{lemme}
This last subsection is devoted to the following proposition.
Theorem \ref{theoreme} follows from this proposition, Proposition
\ref{proposition1} and Proposition \ref{proposition2}.
\begin{prop}
Let us denote
\[
\tilde{A}(t,x)=\frac{1}{2\pi}\int_0^{2\pi}A(R_\theta
x,t)\,d\theta.
\]
Then $a(e)$ is non negative and equal to:
\[
\frac{1}{2\pi n^2}\int_0^\infty(-\partial^2_{tt}\tilde{A})(\frac{-s}{2\pi
n },2\sqrt{e}\sqrt{1-\cos s})\,ds.
\]
\end{prop}

\noindent {\bf Proof.} Thanks to Lemma \ref{equationlimitegyromoyenee} and lemma \ref{hypotheseA}, we have
\[
a(e)=\int_{s=0}^{\infty}\int_{0}^{2\pi}R_{\theta}v\cdot(-\nabla^2_{xx}A)(-\frac{s}{2\pi
n },R_\theta v^\perp-R_{-s}R_\theta v^\perp)\cdot R_{-s}R_\theta v\,ds\,d\theta.
\]
Since
\begin{eqnarray*}
&&\qquad \qquad -\nabla^2_{xx}A(-\frac{s}{2\pi n},v^\perp-R_{-s}
v^\perp)\cdot
R_{-s} v\\
&&=\frac{1}{2\pi n}\nabla_x\partial_sA(-\frac{s}{2\pi n
},v^\perp-R_{-s} v^\perp)+\partial_s(\nabla_xA(-\frac{s}{2\pi n
},v^\perp-R_{-s} v^\perp)),
\end{eqnarray*}
we find
\begin{eqnarray*}
a(e)&=&\frac{1}{2\pi n}\int_{s=0}^\infty\int_{0}^{2\pi}R_\theta v\cdot\nabla_{x}\partial_sA(-\frac{s}{2\pi n},R_\theta v^\perp-R_{\theta-s}v^\perp)\,ds\,d\theta\\
&&-\int_{0}^{2\pi}R_\theta v\cdot\nabla_{x}A(0,0)\,d\theta\\
&=&\frac{1}{2\pi n
}\int_{s=0}^\infty\int_{0}^{2\pi}R_\theta\vec{e}\cdot\nabla_{x}\partial_sA(-\frac{s}{2\pi
n },R_{\theta+\pi/2}\vec{e}-R_{-s+\theta+\pi/2}
\vec{e})\,ds\,d\theta
\end{eqnarray*}
where $\vec{e}=(\sqrt{2e},0)$. Let us do the change of variables
$s'=\theta-s$ to get
\[
a(e)=\frac{1}{2\pi n}\int_0^{2\pi}\int_{\R}{\bf 1}_{\{
s\leq \theta\}}R_\theta\vec{e}\cdot\nabla_x\partial_s
A(\frac{s-\theta}{2\pi n
},R_{\theta+\pi/2}\vec{e}-R_{s+\pi/2}\vec{e})\,d\theta\,ds.
\]
Next we have
\begin{eqnarray*}
R_\theta\vec{e}\cdot\nabla_x\partial_sA(\frac{s-\theta}{2\pi n},R_{\theta+\pi/2}\vec{e}-R_{s+\pi/2}\vec{e})
&=&-\frac{1}{2\pi n}\partial^2_{ss}A(\frac{s-\theta}{2\pi n},R_{\theta+\pi/2}\vec{e}-R_{s+\pi/2}\vec{e})\\
&&-\partial_\theta\left\{\partial_sA(\frac{s-\theta}{2\pi
n},R_{\theta+\pi/2}\vec{e}-R_{s+\pi/2}\vec{e}) \right\}.
\end{eqnarray*}
Integrating by parts the second term of the RHS gives \be
&&\frac{1}{2\pi n}\int_0^{2\pi}\int_{\R} {\bf 1}_{\{ s\leq
\theta\}}\partial_\theta\left\{\partial_sA(\frac{s-\theta}{2\pi
n},R_{\theta+\pi/2}\vec{e}-R_{s+\pi/2}\vec{e})
\right\}\,ds\,d\theta\\
&=&\frac{1}{2\pi n}\int_0^{2\pi}\partial_{\theta}\int_{\R}
{\bf 1}_{\{ s\leq
\theta\}}\partial_sA(\frac{s-\theta}{2\pi n},R_{\theta+\pi/2}\vec{e}-R_{s+\pi/2}\vec{e})\,ds\,d\theta\\
&-&\frac{1}{2\pi n}\int_0^{2\pi}\int_{\R} \delta_{s=\theta}\partial_sA(\frac{s-\theta}{2\pi n},R_{\theta+\pi/2}\vec{e}-R_{s+\pi/2}\vec{e})\,ds\,d\theta\\
&=&\frac{1}{2\pi n}\int_{-\infty}^{2\pi} \partial_sA(\frac{s-2\pi}{2\pi n},R_{\pi/2}\vec{e}-R_{s+\pi/2}\vec{e})\,ds\\
&-&\frac{1}{2\pi n}\int_{-\infty}^0 \partial_sA(\frac{s}{2\pi n},R_{\pi/2}\vec{e}-R_{s+\pi/2}\vec{e})\,ds\\
&-&\frac{1}{2\pi n}\int_0^{2\pi}\partial_s A(0,0)\,ds \en The
first two lines cancel by doing the change of variables
$s'=s-2\pi$ and the third line vanishes thanks to Lemma
\ref{hypotheseA}, thus
\begin{eqnarray*}
a(e)=\frac{1}{(2\pi n)^2}\int_0^{2\pi}\int_{\R} {\bf
1}_{\{ s\leq \theta\}}(-\partial_{ss}^2 A)(\frac{s-\theta}{2\pi
n},R_{\theta+\pi/2}\vec{e}-R_{s+\pi/2}\vec{e})\,d\theta\,ds
\end{eqnarray*}
Doing the change of variables $s'=\theta-s$ gives
\begin{eqnarray*}
a(e)&=&\frac{1}{(2\pi n)^2}\int_0^{2\pi}\int_0^\infty 
(-\partial_{ss}^2 A)(-\frac{s}{2\pi n},R_{\theta+\pi/2}((I-R_{-s})\vec{e}))\,d\theta\,ds\\
&=&\frac{1}{2\pi
n^2}\int_0^\infty (-\partial_{ss}^2\tilde{A})(-\frac{s}{2\pi n},|(I-R_{-s})\vec{e}|)\,ds.
\end{eqnarray*}
Finally
\begin{eqnarray*}
|(I-R_s)\vec{e}|&=&\sqrt{(|1-\cos s|^2+\sin^2 s)}\sqrt{2e}\\
&=&\sqrt{2(1-\cos s)}\sqrt{2e}\\
&=&2\sqrt{e}\sqrt{1-\cos s}\\
&=&2\sqrt{2e}|\sin(s/2)|
\end{eqnarray*}
which ends the proof of the second assertion. $\hfill \Box$
\\
\\
{\bf Computation of the sign of the diffusion coefficient.}
\\
Here we check the non-negativity of the diffusion coefficient by expressing it in another form.
Thanks to lemma \ref{equationlimitegyromoyenee} and to hypothesis $(H3), (H4)$ we have
\be
a(e)=\frac{1}{2N}\int_{s=0}^{+\infty}\int_{-2\pi N}^{2\pi N} R_{\theta} v\cdot
(-\nabla^2_{xx}A) (-\frac{s}{2\pi n},R_{\theta} v^\perp - R_{\theta-s} v^\perp)\cdot R_{\theta-s} v \,ds\, d\theta 
\en
for all $N$. Then doing the change of variable $s:=\theta -s$ we find:
\be
a(e)=\frac{1}{2N}\int_{s\in \R}\int_{-2\pi N}^{2\pi N} {\bf 1}_{\{\theta \geq s\}} R_{\theta} v\cdot(-\nabla^2_{xx}A)(\frac{s-\theta}{2\pi n},R_{\theta} v^\perp - R_{s} v^\perp)\cdot R_s v ds d\theta  
\en
But we remind that thanks to hypothesis $(H4)$ 
\be
-\nabla^2_{xx}A(\frac{s-\theta}{2\pi n},R_{\theta} v^\perp - R_{s} v^\perp)=\lim_{\epsilon\to 0}{\bf E}\left(\nabla_x V^\epsilon(-\frac{s}{2\pi n},- R_{s} v^\perp)\otimes \nabla_x V^\epsilon(-\frac{\theta}{2\pi n},-R_{\theta} v^\perp) \right)
\en
Thus
\be
&&a(e)=\lim_{\epsilon,N}\frac{1}{2N}\cdot \\
&&\int\int_{-2\pi N}^{2\pi N}{\bf 1}_{\{\theta \geq s\}}
{\bf E}\left([\nabla_x V^\epsilon(-\frac{s}{2\pi n},- R_{s} v^\perp)\cdot R_s v][ \nabla_x V^\epsilon(-\frac{\theta}{2\pi n},-R_{\theta} v^\perp)\cdot R_{\theta}v] \right)dsd\theta
\en
Interverting $s$ and $\theta$ we see that we can replace ${\bf 1}_{\{\theta -s\geq 0\}}$ by ${\bf 1}_{\{s-\theta\geq 0\}}$
and finally by adding both we obtain:
\be
a(e)=\lim_{\epsilon, N}\frac{1}{4N}{\bf E}\left(\left[\int_{-2\pi N}^{2\pi N}\nabla_x V^\epsilon{}(-\frac{s}{2\pi n},- R_{s} v^\perp)\cdot R_s v ds\right]^2\right)
\en
which is a positive quantity. $\hfill \Box$

\bigskip
\noindent
{\bf Acknowledgments:} The idea of this paper was born as the authors attended a course
on plasma's turbulence given by Xavier Garbet at the CEA center of of Cadarache. They are glad
to thank him for having introduced them  to the subject and for subsequent 
discussions.


\bibliography{loeper-vasseur-biblio}

\end{document}